\documentclass[review]{elsarticle}

\usepackage{lineno,hyperref}
\modulolinenumbers[1]

\journal{arXiv}

 \biboptions{numbers,sort&compress}
  
\usepackage{graphicx}
\usepackage{amsmath}
\usepackage{amsfonts}
\usepackage[compact]{titlesec}
\usepackage{amsthm}
\usepackage{amssymb}
\usepackage{multirow}
\usepackage{color}
\usepackage{mathtools}
\usepackage{geometry}
\usepackage{times}
\usepackage{algorithm}
\usepackage{algorithmicx, algpseudocode}

\usepackage{verbatim}

\newtheorem{theorem}{Theorem}[section]

\theoremstyle{definition}

\newtheorem{remark}[theorem]{Remark}
\newcommand{\bfs}[1]{{\boldsymbol #1}}

\definecolor{Green}{rgb}{0,.5,0}
\definecolor{Blue}{rgb}{0,.1,.85}
\definecolor{Cyan}{rgb}{.2,.6,.7}
\definecolor{Purple}{rgb}{.5,0,1}
\definecolor{deepred}{rgb}{.8,.1,.2}

\begin{document}

\begin{frontmatter}

\title{Isogeometric spectral approximation for elliptic differential operators} 


\author[ad,ad1]{Quanling Deng\corref{corr}}
\cortext[corr]{Corresponding author}
\ead{Quanling.Deng@curtin.edu.au}

\author[ad,ad1]{Vladimir Puzyrev}
\ead{Vladimir.Puzyrev@curtin.edu.au}

\author[ad,ad1,ad2]{Victor Calo}
\ead{Victor.Calo@curtin.edu.au}

\address[ad]{Curtin Institute for Computation, Curtin University, Kent Street, Bentley, Perth, WA 6102, Australia}

\address[ad1]{Department of Applied Geology, Curtin University, Kent Street, Bentley, Perth, WA 6102, Australia}
\address[ad2]{Mineral Resources, Commonwealth Scientific and Industrial Research Organisation (CSIRO), Kensington, Perth, WA 6152, Australia}

\begin{abstract}
We study the spectral approximation of a second-order elliptic differential eigenvalue problem that arises from structural vibration problems using isogeometric analysis. In this paper, we generalize recent work in this direction. We present optimally blended quadrature rules for the isogeometric spectral approximation of a diffusion-reaction operator with both Dirichlet and Neumann boundary conditions. The blended rules improve the accuracy and the robustness of the isogeometric approximation. In particular, the optimal blending rules minimize the dispersion error and lead to two extra orders of super-convergence in the eigenvalue error. Various numerical examples (including the Schr$\ddot{\text{o}}$dinger operator for quantum mechanics)  in one and three spatial dimensions demonstrate the performance of the blended rules.
\end{abstract}



\end{frontmatter}

\linenumbers

\section{Introduction}  \label{sec:intro}
Differential eigenvalue problems arise in a wide range of applications, such as the vibration of elastic bodies in structural mechanics and the multi-group diffusion in nuclear reactors \cite{strang1973analysis}. In general, analytical solutions for these problems are impossible and numerical methods are used. Numerical methods for approximating these differential eigenvalue problems lead to a generalized matrix eigenvalue problem, which is then solved numerically. Different numerical methods result in different matrices and the widely-used methods include finite elements \cite{strang1973analysis,bramble1973rate,osborn1975spectral, chatelin1983spectral, babuvska1991eigenvalue, canuto1978eigenvalue,mercier1978eigenvalue,mercier1981eigenvalue}, isogeometric elements \cite{hughes2014finite, calo2017dispersion, puzyrev2017dispersion, deng2017dispersion,deng2018dispersion,bartovn2017generalization,calo2017quadrature,puzyrev2017spectral}, discontinuous Galerkin (DG) \cite{antonietti2006discontinuous,giani2015hp}, hybridizable discontinuous Galerkin (HDG) \cite{gopalakrishnan2015spectral}, and a recently developed hybrid high-order (HHO) method \cite{calo2017spectral}.

Early work \cite{strang1973analysis,bramble1973rate,osborn1975spectral} used conforming finite elements on simplicial meshes and the method demonstrated convergence rates of order $h^{2p}$ for the eigenvalues and of order $h^p$ in the energy norm for the eigenfunctions provided that the eigenfunctions are smooth enough. Sharp and optimal estimates of the numerical eigenfunctions and eigenvalues of finite element analysis are established in \cite{banerjee1989estimation,banerjee1992note,banerjee1992analysis}.  In \cite{canuto1978eigenvalue,mercier1978eigenvalue,mercier1981eigenvalue}, the authors studied the spectral approximation of elliptic operators by mixed and mixed-hybrid methods and optimal error estimates were established. Similar results were obtained more recently in \cite{antonietti2006discontinuous,giani2015hp,gopalakrishnan2015spectral,calo2017spectral} using the DG, HDG, and HHO methods. The spectral approximation by the HDG method leads to a convergence of order $h^{2p+1}$ for the eigenvalues.  A non-trivial post-processing, which utilizes a Rayleigh quotient, is also examined in \cite{gopalakrishnan2015spectral} numerically which leads to an improved convergence of order $h^{2p+2}$ for $p\ge1$. The HHO approximation delivers a convergence of order $h^{2p+2}$ for the eigenvalue errors for all polynomial degrees ($p\ge0$).

Isogeometric analysis is a numerical method introduced in 2005 \cite{hughes2005isogeometric,cottrell2009isogeometric}. The spectral approximation of the elliptic operators arising in structural vibrations were investigated using isogeometric analysis in \cite{cottrell2006isogeometric,hughes2008duality} and the method shows improved spectral approximations over the classical finite elements \cite{cottrell2009isogeometric}. In \cite{hughes2008duality}, a duality principle, which induces a bijective map from spectral analysis to dispersion analysis, was established, which unifies the spectral analysis for structural vibrations (eigenvalue problems) and the dispersion analysis for wave propagations. Further advantages of the method on spectral approximation properties are investigated in \cite{hughes2014finite}.

The recent work in \cite{calo2017dispersion, puzyrev2017dispersion} studies both theoretically and numerically the optimally blended quadrature rules \cite{ainsworth2010optimally} for the isogeometric analysis of the Laplace eigenvalue problem.  In \cite{calo2017dispersion}, the authors establish for $p=1,2,\cdots,7$ the super-convergence of order $h^{2p+2}$ for the eigenvalue errors while maintaining optimal convergence of orders $h^p$ and $h^{p+1}$ for the eigenfunction errors in the $H^1$-seminorm and in the $L^2$-norm, respectively. The work \cite{deng2017dispersion} introduces the dispersion-minimized mass for isogeometric analysis and generalizes the results to arbitrary polynomial degree $p$. In \cite{puzyrev2017spectral}, the authors study the optimally blended quadratures for isogeometric analysis with variable continuity. To reduce the computational costs, \cite{deng2018dispersion} describes new quadrature rules to replace the optimal blending rules. For the source problems, optimal (Gaussian) quadrature rules were proposed for isogeometric analysis in \cite{bartovn2016optimal,bartovn2016gaussian,bartovn2017gauss}.  

In this work, we generalize the work in \cite{calo2017dispersion, puzyrev2017dispersion} to include the reaction effects in the differential operator as well as to consider different boundary conditions. We study numerically the optimal blending quadratures for the generalized differential operator with both Dirichlet and Neumann boundary conditions. We apply the blending rules to approximate the spectrum of the Schr$\ddot{\text{o}}$dinger operator.

The outline of the rest of this paper is as follows. We first describe the model problem and the isogeometric discretization in Section \ref{sec:ps}. We introduce classical and blended quadrature rules in Section \ref{sec:blending}. A brief dispersion error estimations is presented in this section. 
Numerical examples are given in Section \ref{sec:num}. Finally, Section \ref{sec:con} summarizes our findings and describes future research directions.

\section{Problem statement} \label{sec:ps}
We consider the second-order differential eigenvalue problem: Find the eigenpair $(\lambda, u)$ such that
\begin{equation} \label{eq:pde}
\begin{aligned}
- \Delta u + \gamma u & = \lambda u \quad  \text{in} \quad \Omega, \\
u & = 0 \quad \text{on} \quad \partial \Omega,
\end{aligned}
\end{equation}
where $\Delta = \nabla^2$ is the Laplacian, $\gamma = \gamma(x) \in L^2(\Omega)$ is a smooth and non-negative function, and $\Omega \subset  \mathbb{R}^d, d=1,2,3$ is a bounded open domain with Lipschitz boundary. This problem is a Sturm-Liouville eigenvalue problem (see, for example, \cite{strang1973analysis,evans2010partial}) which has a countable infinite set of eigenvalues ${\lambda _j} \in {\mathbb{R}}$ 
\begin{equation}
0 < {\lambda _1} < {\lambda _2} \le \cdots \le {\lambda _j} \le \cdots
\end{equation}
with an associated set of orthonormal eigenfunctions ${u_j}$ 
\begin{equation}
({u_j},{u_k}) = \int_{\Omega}  {{u_j}(x){u_k}} (x) \ \text{d} \bfs{x} = {\delta _{jk}},
\end{equation}
where $\delta _{jk}$ is the Kronecker delta which is equal to 1 when $j=k$ and 0 otherwise. The set of all the eigenvalues is the spectrum of the operator. We normalize the eigenfunctions in the $L^2$ space and hence the eigenfunctions are orthonormal with each other under the scalar inner product. Now, let us define two bilinear forms 
\begin{equation} \label{eq:bfs}
a(w,v) = \int_{\Omega}  \nabla w \cdot \nabla v + \gamma w v \ \text{d} \bfs{x} \quad \text{and} \quad b(w,v) = (w, v) = \int_{\Omega} w v \ \text{d} \bfs{x}, \quad \forall w, v \in H^1_0(\Omega),
\end{equation}
where $H^1_0(\Omega)$ is the Sobolev space with functions vanishing at the boundary $\partial \Omega.$ 
These two inner products are associated with the following energy and $L^2$ norms
\begin{equation}
\left\| w \right\|_E = \sqrt{a(w,w)}, \qquad \left\| w \right\| = \| w \|_{L^2(\Omega)} = \sqrt{(w,w)}.
\end{equation}
Using this notation, the eigenfunctions are also orthogonal with each other under the energy inner product, that is, 
\begin{equation} \label{eq:oe}
a(u_j, u_k) = (-\Delta {u_j} + \gamma u_j, {u_k}) = ({\lambda _j}{u_j},{u_k}) = {\lambda _j} ({u_j},{u_k})  = {\lambda _j}{\delta _{jk}},
\end{equation}
where we have used the integration by parts on \eqref{eq:pde} and selected the weighting functions to be  eigenfunctions.

At the continuous level, the weak formulation for the eigenvalue problem \eqref{eq:pde} is: Find all eigenvalues $\lambda \in {\mathbb{R}}$ and eigenfunctions $u \in V$ such that,
\begin{equation} \label{eq:weak}
a(w, u) = \lambda b(w, u), \quad \forall \ w \in H^1_0(\Omega),
\end{equation}
while at the discrete level, the isogeometric analysis for the eigenvalue problem \eqref{eq:pde} is: Find all eigenvalues $\lambda_h \in {\mathbb{R}}$ and eigenfunctions $u_h \in V_h$ such that,
\begin{equation} \label{eq:vf}
a(w_h, u_h) = \lambda_h b(w_h, u_h), \quad \forall \ w_h \in V_h(\Omega),
\end{equation}
where $V_h \subset H^1_0(\Omega)$ is the solution and test space, which is spanned by the B-spline or non-uniform rational basis spline (NURBS) basis functions.

Following \cite{de1978practical,piegl2012nurbs}, the definition of the B-spline basis functions in one dimension is as follows. 
Let $X = \{x_0, x_1, \cdots, x_m \}$ be a knot vector with knots $x_j$, that is, a nondecreasing sequence of real numbers called knots.  The $j$-th B-spline basis function of degree $p$, denoted as $\phi^j_p(x)$, is defined as 
\begin{equation} \label{eq:B-spline}
\begin{aligned}
\phi^j_0(x) & = 
\begin{cases}
1, \quad \text{if} \ x_j \le x < x_{j+1} \\
0, \quad \text{otherwise} \\
\end{cases} \\ 
\phi^j_p(x) & = \frac{x - x_j}{x_{j+p} - x_j} \phi^j_{p-1}(x) + \frac{x_{j+p+1} - x}{x_{j+p+1} - x_{j+1}} \phi^{j+1}_{p-1}(x).
\end{aligned}
\end{equation}

In this paper, we use the B-splines on uniform meshes with non-repeating knots, that is, we use B-splines with maximum continuity.
We approximate the eigenfunction as a linear combination of the B-spline basis functions.
Using linearity and substituting all the B-spline basis functions for $w_h$ in \eqref{eq:vf} leads to the matrix eigenvalue problem
\begin{equation} \label{eq:mevp}
\mathbf{K} \mathbf{U} = \lambda^h \mathbf{M} \mathbf{U},
\end{equation}
where $\mathbf{K}_{ab} =  a(\phi_a, \phi_b), \mathbf{M}_{ab} = b(\phi_a, \phi_b),$ and $\mathbf{U}$ is the corresponding representation of the eigenvector as the coefficients of the B-spline basis functions. For simplicity, the matrix $\mathbf{K}$ (although it contains a scaled mass) is referred as the stiffness matrix while the matrix $\mathbf{M}$ is referred as the mass matrix, and $(\lambda_h, \mathbf{u}_h)$ is the unknown eigenpair.

\section{Blending quadratures and dispersion errors} \label{sec:blending}
In this section, we present the quadratures as well as their optimal blendings. Following earlier work \cite{calo2017dispersion, puzyrev2017dispersion} and its recent generalization \cite{deng2017dispersion}, we omit the details to briefly give the dispersion errors for the quadrature rules. The blending rules are optimal in the sense of delivering minimal dispersion error.

\subsection{Quadrature rules}
In practice, we evaluate the integrals involved in $a(w_h, u_j^h) $ and $b(w_h, u_j^h)$ numerically, that is, approximated by quadrature rules. On a reference element $\hat K$, a quadrature rule is of the form
\begin{equation} \label{eq:qr}
\int_{\hat K} \hat f(\hat{\boldsymbol{x}}) \ \text{d} \hat{\boldsymbol{x}} = \sum_{l=1}^{n} \hat{\varpi}_l \hat f (\hat{n_l}) +  \hat{E}_n,
\end{equation}
where $\hat{\varpi}_l$ are the weights, $\hat{n_l}$ are the nodes, $n$ is the number of quadrature points, and $\hat{E}_n$ is the error of the quadrature rule. For each element $K$, we assume that there is an invertible map $\sigma$ such that $K = \sigma(\hat K)$, which determines the correspondence between the functions on $K$ and $\hat K$. Assuming $J_K$ is the Jacobian of the mapping, \eqref{eq:qr} induces a quadrature rule over the element $K$ given by
\begin{equation} \label{eq:q}
\int_{K}  f(\boldsymbol{x}) \ \text{d} \boldsymbol{x} = \sum_{l=1}^{n} \varpi_{l,K} f (n_{l,K}) + E_n,
\end{equation}
where $\varpi_{l,K} = \text{det}(J_K) \hat \varpi_l$ and $n_{l,K} = \sigma(\hat n_l)$.

The quadrature rule is exact for a given function $f(x)$ when the remainder $E_n$ is exactly zero. For simplicity, we denote by $G_m$ the $m-$point Gauss-Legendre quadrature rule,  by $L_m$ the $m-$point Gauss-Lobatto quadrature rule, and by $O_p$ the optimal blending scheme for the $p$-th order isogeometric analysis with maximum continuity. In one dimension, $G_m$ and $L_m$ fully integrate polynomials of order $2m-1$ and $2m-3$, respectively (see, for example, \cite{stoer2013introduction, bartovn2016optimal}).

Applying the quadrature rules to \eqref{eq:vf}, we have the approximated form
\begin{equation} \label{eq:vfho}
 a_h(w_h, \tilde u_h) =  \tilde\lambda_h  b_h(w_h, \tilde u_h), \quad \forall \ w_h \in V_h,
\end{equation}
where for $w,v \in V_h$
\begin{equation} \label{eq:ba}
 a_h(w, v) = \sum_{K \in \mathcal{T}_h} \Big( \sum_{l=1}^{n_1}  \varpi_{l,K}^{(1)} \nabla w (n_{l,K}^{(1)} ) \cdot \nabla v (n_{l,K}^{(1)} ) + \sum_{l=1}^{n_2}  \varpi_{l,K}^{(2)} \gamma (n_{l,K}^{(2)} )  w (n_{l,K}^{(2)} ) v (n_{l,K}^{(2)} ) \Big)
\end{equation}
and
\begin{equation} \label{eq:bb}
 b_h(w, v) = \sum_{K \in \mathcal{T}_h} \sum_{l=1}^{n_3} \varpi_{l,K}^{(3)} w (n_{l,K}^{(3)} ) v (n_{l,K}^{(3)} ),
\end{equation}
where $\{\varpi_{l,K}^{(j)}, n_{l,K}^{(j)} \}$ with $j=1,2,3$  specifies three (possibly different) quadrature rules. Here, we assume that we apply the same quadrature rules for the $L^2$ inner products in \eqref{eq:ba} and \eqref{eq:bb}. With these quadrature rules, we can rewrite (with slight abuse of the notation) the matrix eigenvalue problem \eqref{eq:mevp} as
\begin{equation} \label{eq:amevp}
\mathbf{K} \tilde{\mathbf{U}} = \tilde \lambda^h \mathbf{M} \tilde{\mathbf{U}},
\end{equation}
where $\mathbf{K}_{ab} =  a_h(\phi_a, \phi_b), \mathbf{M}_{ab} = b_h(\phi_a, \phi_b),$ and $\tilde{\mathbf{U}}$ is the corresponding representation of the eigenvector as the coefficients of the basis functions.

\begin{remark}
For multidimensional problems on tensor product grids, the stiffness and mass matrices can be expressed as Kronecker products of 1D matrices \cite{gao2014fast}. For example, in the 2D case, assume that $\gamma$ is a constant. We define $\gamma=\gamma_{2D} = 2\gamma_{1D}$ and let $\phi^i_p(x) \phi^j_p(y)$ and $\phi^k_p(x) \phi^l_p(y)$ be two 2D basis functions. Using the definition \eqref{eq:bfs}, we calculate 
\begin{equation} \label{eq:km}
\begin{aligned}
\mathbf{K} _{i j k l} & = a\big(\phi^i_p(x) \phi^j_p(y), \phi^k_p(x) \phi^l_p(y) \big) \\
& = \int_{\Omega}  \nabla \big(\phi^i_p(x) \phi^j_p(y)\big) \cdot \nabla \big( \phi^k_p(x) \phi^l_p(y) \big) + \gamma_{2D} \phi^i_p(x) \phi^j_p(y) \phi^k_p(x) \phi^l_p(y)  \ \text{d} \bfs{x} \\
& = \int_X \Big( \frac{\partial \phi^i_p(x) }{\partial x} \frac{\partial \phi^k_p(x) }{\partial x} + \gamma_{1D} \phi^i_p(x) \phi^k_p(x) \Big) \ \text{d} x  \int_Y \phi^j_p(y) \phi^l_p(y)  \ \text{d} y\\ 
& \quad + \int_X \phi^i_p(x) \phi^k_p(x)  \ \text{d} x \int_Y \Big( \frac{\partial \phi^j_p(y) }{\partial y} \frac{\partial \phi^l_p(y) }{\partial y} + \gamma_{1D} \phi^j_p(y) \phi^l_p(y)  \Big)\ \text{d} y \\
& = \mathbf{K} ^{1D} _{i k} \mathbf{M} ^{1D} _{j l} + \mathbf{M} ^{1D} _{i k}  \mathbf{K} ^{1D} _{j l},
\end{aligned}
\end{equation}
where $X$ and $Y$ specify the intervals of each dimension in $\Omega$. Similarly, we obtain 
\begin{equation} \label{eq:2d-1d}
\mathbf{M} _{i j k l} = \mathbf{M} ^{1D} _{i k} \mathbf{M} ^{1D} _{j l}.
\end{equation}
Herein, $\mathbf{M} ^{1D} _{i j}$ and $\mathbf{K} ^{1D} _{i j}$ are the mass and stiffness matrices of the 1D problem with $\gamma = \gamma_{1D}$ in \eqref{eq:pde}. We refer the reader to \cite{de2007grid} for the description of the summation rules. 

\end{remark}


\subsection{Blended quadratures}
Given two quadrature rules $Q_1 = \{\varpi_{l,K}^{(1)}, n_{l,K}^{(1)} \}_{l=1}^{n_1}$ and $Q_2 = \{\varpi_{l,K}^{(2)}, n_{l,K}^{(2)} \}_{l=1}^{n_2}$, the blended quadrature rule, denoted as $Q_\tau$, is defined as
\begin{equation}
Q_\tau = \tau Q_1 + (1-\tau) Q_2,
\end{equation}
where $\tau$ is referred as the blending parameter. Applying the blended rule $Q_\tau$ for the integration of a function $f$, we have
\begin{equation}
\begin{aligned}
\int_{K}  f(\boldsymbol{x}) \ \text{d} \boldsymbol{x} & = \tau \Big( \sum_{l=1}^{n_1} \varpi_{l,K}^{(1)} f (n_{l,K}^{(1)}) + E_{n_1}^{Q_1} \Big) + (1 - \tau) \Big( \sum_{l=1}^{n_2} \varpi_{l,K}^{(2)} f (n_{l,K}^{(2)}) + E_{n_2}^{Q_2} \Big) \\
& = \tau \sum_{l=1}^{n_1} \varpi_{l,K}^{(1)} f (n_{l,K}^{(1)}) + (1 - \tau) \sum_{l=1}^{n_2} \varpi_{l,K}^{(2)} f (n_{l,K}^{(2)})  + \Big( \tau E_{n_1}^{Q_1} + (1-\tau) E_{n_2}^{Q_2} \Big). \\
\end{aligned}
\end{equation}
Thus, the error for the blending rule is the same as blending of the errors, that is,
\begin{equation} \label{eq:eQtau}
E^{Q_\tau} = \tau E_{n_1}^{Q_1} + (1-\tau) E_{n_2}^{Q_2}.
\end{equation}
Assuming that $Q_1$ and $Q_2$ integrate polynomials up to order $k_1$ and $k_2$, respectively, \eqref{eq:eQtau} shows that the blending rule integrates polynomials up to order $\min(k_1,k_2)$. For example, in one dimension, the blending rule 
\begin{equation}
Q_\tau = \tau G_m + (1-\tau) L_m
\end{equation}
integrates polynomials up to order $2m-3.$

For the dispersion analysis on the Helmholtz equation ($\gamma=0, \lambda = \omega^2$ in \eqref{eq:pde}), the blended rule shows smaller dispersion errors. In fact, the optimal blending of spectral elements and finite elements, which is realized by optimally blended quadratures, leads to two extra order of super-convergence on the dispersion error; see \cite{ainsworth2010optimally}.  This fact motivates the work (see \cite{calo2017dispersion, puzyrev2017dispersion, deng2017dispersion}) of finding the optimal blending rules for the isogeometric analysis for differential eigenvalue problems \eqref{eq:pde} with $\gamma=0$. In the following section, we present the dispersion error-minimized blending rules.

\subsection{Dispersion errors and optimal blending quadratures}
Following earlier work \cite{calo2017dispersion,deng2017dispersion}, based on the dual principle in \cite{hughes2008duality}, the dispersion errors of the isogeometric elements using quadratures can be characterized by the eigenvalue errors. For simplicity, we assume that $\gamma=0$. For $C^1$ quadratic isogeometric elements (for linear elements, it is the same with the finite element case and we refer the readers to \cite{puzyrev2017dispersion}), the relative errors are
\begin{equation}
\begin{aligned}
\frac{\lambda_h^{G_3} - \lambda}{\lambda}  & = \frac{1}{720} \Lambda^4 + \mathcal{O}(\Lambda^6), \\
\frac{\lambda_h^{L_3} - \lambda}{\lambda}  & = -\frac{1}{1440} \Lambda^4 + \mathcal{O}(\Lambda^6), \\
\end{aligned}
\end{equation}
where $\Lambda = \omega h$ with $\omega^2 = \lambda$ and $\lambda_h^{Q}$ denotes the approximate eigenvalue while using the quadrature rule $Q$. 
The blending of these two rules, that is, $Q_\tau = \tau G_3 + (1-\tau) L_3$,  leads to the error representation 
\begin{equation}
\begin{aligned}
\frac{\lambda_h^{Q_\tau} - \lambda}{\lambda}  & = \frac{2 - 3 \tau }{1440} \Lambda^4 + \mathcal{O}(\Lambda^6). \\
\end{aligned}
\end{equation}
For $\tau = 2/3$, we obtain the two extra orders in the error representation and we call this case the optimal blending. The error representation of the optimal blending is 
\begin{equation}
\begin{aligned}
\frac{\lambda_h^{O_p} - \lambda}{\lambda}  & = \frac{11 }{60480} \Lambda^6 + \mathcal{O}(\Lambda^8). \\
\end{aligned}
\end{equation}

For $C^2$ cubic elements, the optimal blending parameter is $\tau = -3/2$ and we refer to \cite{calo2017dispersion} for $p\le 7$ and \cite{deng2017dispersion} for the general case.
The convergence rate for eigenpairs computed using isogeometric elements is $O\left( \Lambda^{2p} \right) $ as shown in \cite{cottrell2006isogeometric}. The optimal blending leads to a $O\left( \Lambda^{2p+2} \right) $ convergence rate for the relative eigenvalue errors.

For a constant $\gamma$, we redefine the eigenvalue problem \eqref{eq:pde} as
\begin{equation} \label{eq:pdee}
\begin{aligned}
- \Delta u & = \hat \lambda u \quad  \text{in} \quad \Omega, \\
u & = 0 \quad \text{on} \quad \partial \Omega,
\end{aligned}
\end{equation}
where $\hat \lambda = \lambda - \gamma$. Once the eigenvalue problem \eqref{eq:pdee} is solved using isogeometric elements with optimal blending rules, we post-process the approximated eigenvalue of \eqref{eq:pde} as 
\begin{equation}
\lambda_h = \hat \lambda_h + \gamma.
\end{equation}

In Section \ref{sec:num}, we present the numerical studies of the equation \eqref{eq:pde} for variable coefficient $\gamma$ using isogeometric analysis with the optimal blending rules.

\section{Numerical examples} \label{sec:num}

In this section, we present numerical examples of the one and three dimensional problems described in Section 2 to show how optimal quadratures reduce the approximation errors in isogeometric analysis.

\subsection{1D results}
The 1D elliptic eigenvalue problem \eqref{eq:pde} with $\gamma=0$ and homogeneous Neumann boundary conditions has the exact eigenpairs $\lambda _j = j^2 \pi ^2, u_j = \sqrt 2 \cos (j\pi x), j = 1, 2, \cdots.$ The approximate eigenvalues $\lambda _j^h$ are sorted in ascending order and are compared to the corresponding exact eigenvalues ${\lambda _j}$.

\begin{figure}[!ht]
\centering\includegraphics[width=1.0\linewidth]{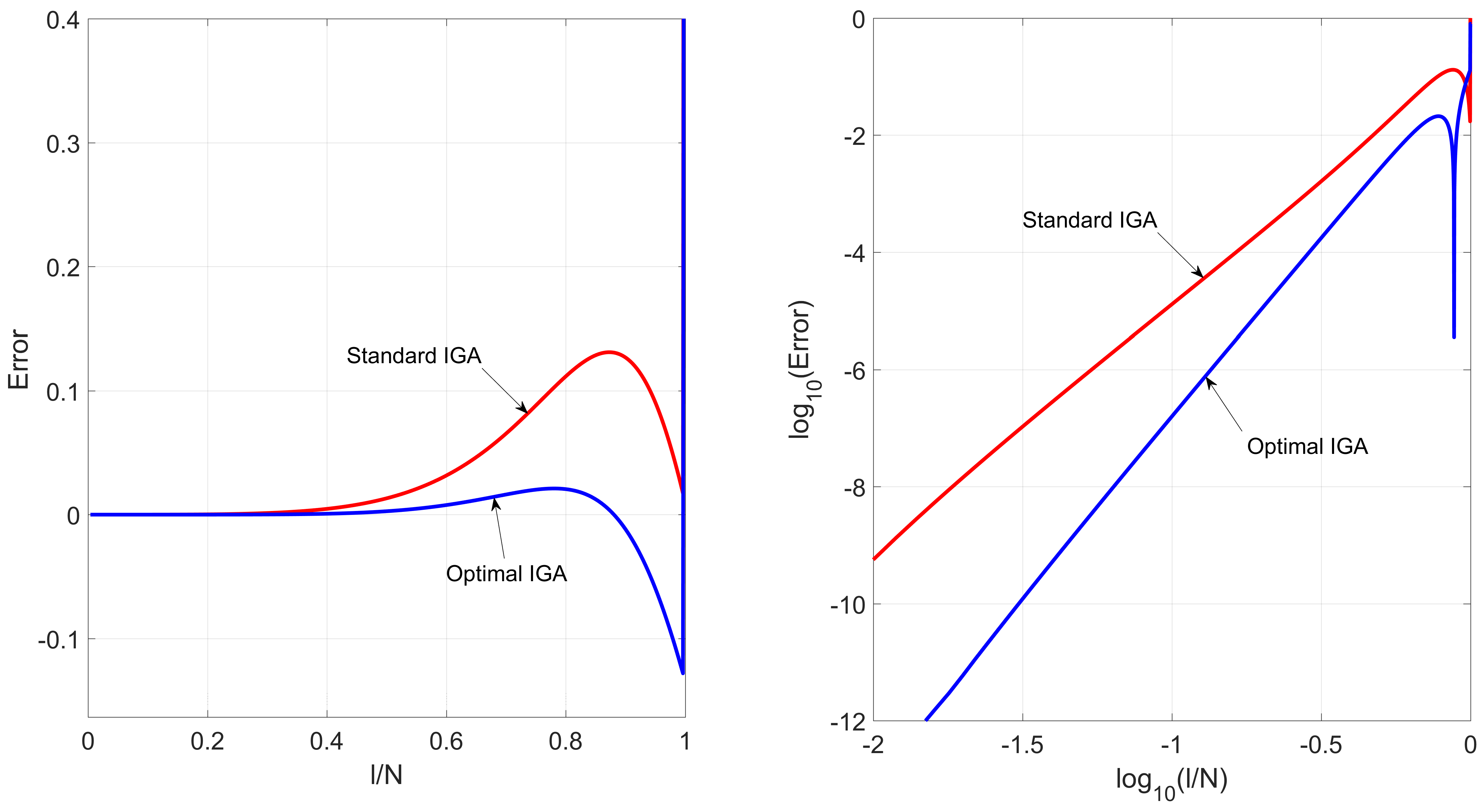}
\caption{Approximation errors for $C^1$ quadratic isogeometric elements with standard Gauss quadrature rule and optimal rule on linear (left) and logarithmic scales (right). The total number of degrees of freedom (discrete modes) is $N = 1000$.}
\label{fig:1}
\end{figure}

Figure \ref{fig:1} compares the approximation errors of $C^1$ quadratic isogeometric elements using the standard Gaussian quadrature and the optimal rule for problem \eqref{eq:pde} with homogeneous Neumann boundary conditions. The use of the optimal quadrature leads to more accurate results. The optimal ratio of blending of the Lobatto and Gauss quadrature rules in this case is 2:1 ($\tau=2/3$), which in this particular case coincides with the ratio proposed by Ainsworth and Wajid \cite{ainsworth2010optimally} for finite-spectral elements of the same polynomial order. This ratio is different for higher order isogeometric elements \citep{calo2017dispersion}.

\begin{figure}[!ht]
\centering\includegraphics[width=1.0\linewidth]{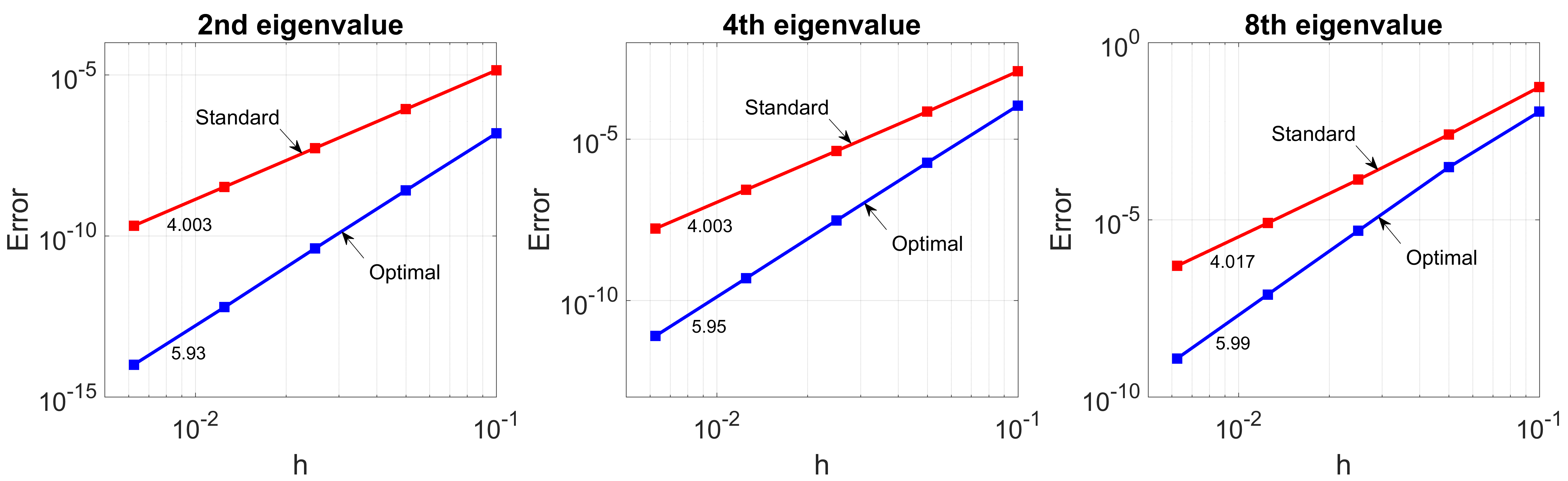}
\caption{Convergence of the errors in the eigenvalue approximation using $C^1$ quadratic isogeometric elements with standard and optimal quadratures. The second (left), fourth (middle) and eighth (right) eigenvalues are shown.}
\label{fig:2}
\end{figure}

Figure \ref{fig:2} shows the convergence of the errors in the eigenvalue approximation with $C^1$ quadratic isogeometric elements. The optimal quadrature rule has two extra orders of convergence in the eigenvalue errors compared to the standard fully-integrated isogeometric elements. Not only the convergence rate, but also the errors are significantly lower for the optimal rule.


\subsection{3D results}
Next, we continue our study with the dispersion properties of the three-dimensional eigenvalue problem \eqref{eq:pde} on tensor product meshes. Optimal methods for multidimensional problems  with constant coefficients and affine mappings can be formed by tensor product of the 1D mass and stiffness matrices \eqref{eq:2d-1d}. The exact eigenvalues and eigenfunctions of the 3D eigenvalue problem are given by
\begin{equation}
{{\lambda _{klm}} = {(k^2+l^2+m^2)}{\pi ^2},\ \ \ {u_{klm}} = 2 \sin (k\pi x) \sin (l\pi y) \sin (m\pi z),}
\end{equation}
for {$k,l,m = 1, 2, ...$}.

\begin{figure}[!ht]
\centering\includegraphics[width=1.0\linewidth]{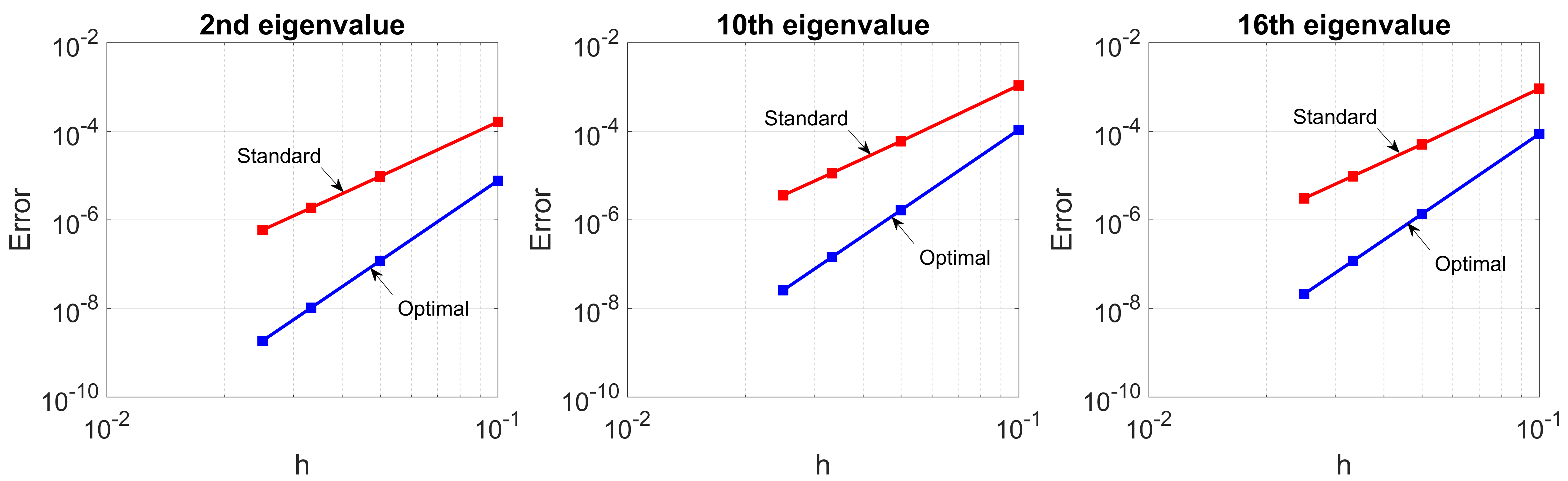}
\caption{Convergence of the errors in the eigenvalue approximation using $C^1$ quadratic isogeometric elements with standard and optimal quadratures. The second (left), tenth (middle) and sixteenth (right) eigenvalues are shown.}
\label{fig:3}
\end{figure}

Figure \ref{fig:3} shows the dispersion errors in the eigenvalue approximation with $C^1$ quadratic  isogeometric elements. Similar to the 1D case, the optimal scheme has two extra orders of convergence in the eigenvalue errors.

\begin{figure}[!ht]
\centering\includegraphics[width=1.0\linewidth]{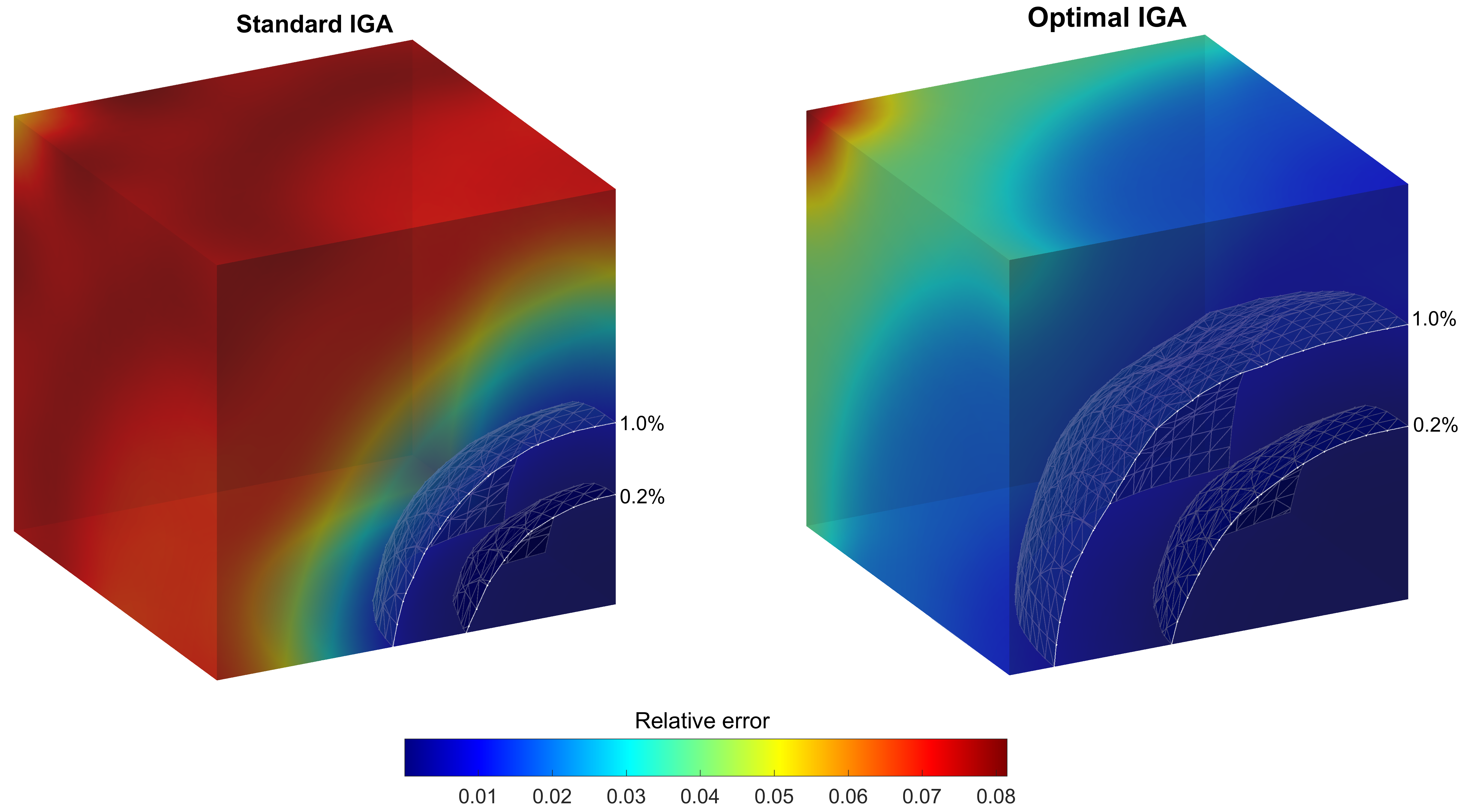}
\caption{Approximation errors for $C^1$ quadratic isogeometric elements with standard Gauss (left) and optimal quadrature rule (right). Color represents the absolute value of the relative error. Isosurfaces show 0.2\% and 1.0\% levels of the relative error.}
\label{fig:4}
\end{figure}

Figure \ref{fig:4} compares the eigenvalue errors of the standard Gauss rule using $C^1$ quadratic elements with the optimal scheme ($\tau=2/3$). The latter has significantly better approximation properties in the entire domain.

These results demonstrate that the use of optimal quadratures in isogeometric analysis significantly improves the accuracy of the discrete approximations compared to the fully-integrated Gauss-based method.

\subsection{Spectral approximation of Schr$\ddot{\text{o}}$dinger operator}
Following the analytical work on Schr$\ddot{\text{o}}$dinger operators in \cite{ciftci2005construction}, we study their numerical approximations in this subsection. We consider the 1D Schr$\ddot{\text{o}}$dinger equation of a quantum particle trapped by the P$\ddot{\text{o}}$schl-Teller potential 
\begin{equation} \label{eq:ptp}
\kappa^2 \Big( \frac{\alpha(\alpha +1)}{\cos^2(\kappa y)} + \frac{\beta(\beta +1)}{\sin^2(\kappa y)} \Big), \qquad 0 < \kappa y < \frac{\pi}{2}, \ \alpha, \beta > 0.
\end{equation}
Applying the scaling $x = \kappa y$, the eigenvalue problem reads: Find the eigenpair $(\lambda, u)$ such that 
\begin{equation}
\begin{aligned}
- \frac{d^2 u}{d x^2} + \Big( \frac{\alpha(\alpha +1)}{\cos^2 x} + \frac{\beta(\beta +1)}{\sin^2 x} \Big) u & = \lambda u,  \qquad 0 < x < \frac{\pi}{2}, \\
u(0) & = 0, \\
u(\frac{\pi}{2}) & = 0,
\end{aligned}
\end{equation}
where we choose $\alpha = \beta = 1$ for simplicity. This eigenvalue problem has the true eigenvalues (see for example \cite{ciftci2005construction})
\begin{equation}
\lambda = (4 + 2j)^2, \quad j=0,1,2, \cdots.
\end{equation}
\begin{table}[h!]
\centering 
\begin{tabular}{| c | c | c | c | c | c | c | c| c| c| c| c|}
\hline
 \multicolumn{2}{|c|}{Set} &  \multicolumn{2}{c|} {$ |\lambda_1^h - \lambda_1| / \lambda_1$ }  & \multicolumn{2}{c|} {$ |\lambda_2^h - \lambda_2| /  \lambda_2$}  & \multicolumn{2}{c|} {$ |\lambda_4^h - \lambda_4| /  \lambda_4$}   \\[0.1cm] \hline
$p$ &  $N$ & $G_{p+1}$ & $O_p$ & $G_{p+1}$ & $O_p$ & $G_{p+1}$ & $O_p$  \\[0.1cm] \hline
& 40	& 3.19e-3	& 6.60e-4	& 1.06e-2	& 1.65e-3	& 3.95e-2	& 3.81e-3 \\[0.1cm]
1& 80	& 7.41e-4	& 8.43e-5	& 2.49e-3	& 2.19e-4	& 9.33e-3	& 5.97e-4 \\[0.1cm]
& 160	& 1.78e-4	& 1.06e-5	& 6.04e-4	& 2.80e-5	& 2.27e-3	& 8.07e-5 \\[0.1cm] \hline
\multicolumn{2}{|c|}{$\rho_1$}	& 2.08	& 2.98	& 2.07	& 2.94	& 2.06	& 2.78 \\[0.1cm] \hline

& 10	& 1.63e-3	& 2.65e-4	& 1.68e-2	& 4.29e-3	& 1.02e+0	& 2.73e-1 \\[0.1cm]
2& 20	& 7.94e-5	& 2.39e-6	& 6.68e-4	& 6.54e-5	& 9.07e-3	& 1.95e-3 \\[0.1cm]
& 40	& 4.62e-6	& 1.11e-7	& 3.61e-5	& 5.24e-7	& 4.07e-4	& 2.83e-5 \\[0.1cm] \hline
\multicolumn{2}{|c|}{$\rho_2$}	& 4.23	& 5.61	& 4.43	& 6.50	& 5.64	& 6.62 \\[0.1cm] \hline

 \end{tabular}
\caption{Relative eigenvalue (EV) errors for $C^0$ linear and $C^1$ quadratic isogeometric elements with Gauss rule $G_{p+1}$ and optimally-blended rules $Q_p$.}
\label{tab:sch} 
\end{table}

Table \ref{tab:sch} shows the relative eigenvalue errors for the first, second and fourth eigenmodes. We present the errors while using both the Gauss rule and optimally blended rule. Here, since the P$\ddot{\text{o}}$schl-Teller potential blows up at the points $x=0,\frac{\pi}{2}$ and the Lobatto rules utilize the interval element end knots as quadrature points, we use the $G_{p+1}$ and $G_p$ optimally blended rules (alternatively, one can use the equivalent nonstandard quadratures; see \cite{calo2017dispersion,deng2018dispersion} for details). The table shows that the eigenvalue errors converge in an order of $2p$ when using the $(p+1)$-point Gauss rule while the error converges in an order of $2p+1$ and $2p+2$ when using the optimal rule $O_p$ for $p=1$ and $p=2$, respectively. The optimal rules were developed for operators with constant coefficients. It is still an open question to develop optimal rules for the operators with variable coefficients. Herein, for the Schr$\ddot{\text{o}}$dinger operator with variable P$\ddot{\text{o}}$schl-Teller potential, the optimal rules improve the eigenvalue errors significantly but the two-extra orders of convergence are not ensured.



%

\section{Conclusions and future outlook} \label{sec:con}
We apply the optimally-blended quadrature rules to approximate the spectrum of a general elliptic differential operator where we account for reaction effects. We show that the optimally blended rules lead to two extra orders of convergence in the eigenvalue errors for both 1D and 3D examples. 

One future direction is the study on the non-uniform meshes and non-constant coefficient differential eigenvalue problems. The study with variable continuity of the B-spline basis functions is also of interest. We will study the dispersion properties of variable continuity in the basis functions on isogeometric elements and study how the dispersion can be minimized by designing goal-oriented quadrature rules.

\section*{Acknowledgement}
This publication was made possible in part by the CSIRO Professorial Chair in Computational Geoscience at Curtin University and the Deep Earth Imaging Enterprise Future Science Platforms of the Commonwealth Scientific Industrial Research Organisation, CSIRO, of Australia. Additional support was provided by the European Union's Horizon 2020 Research and Innovation Program of the Marie Sk{\l}odowska-Curie grant agreement No. 644202, the Mega-grant of the Russian Federation Government (N 14.Y26.31.0013) and the Curtin Institute for Computation. The J. Tinsley Oden Faculty Fellowship Research Program at the Institute for Computational Engineering and Sciences (ICES) of the University of Texas at Austin has partially supported the visits of VMC to ICES. 


\bibliography{igaref}


\end{document}